\documentclass[11pt]{article}
\usepackage{amssymb,amsmath,euscript}
\newtheorem{proposition}{Proposition}[section]

\newtheorem{lemma}[proposition]{Lemma}
\newtheorem{theorem}[proposition]{Theorem}

\newtheorem{corollary}[proposition]{Corollary}

\def\I{{\rm 1 \hskip-2.9truept l}}

\def\l{{\langle}}
\def\r{\rangle}
\def\dim{{\rm dim}_{_{\rm H}}}

\def\dimp{{\rm dim}_{_{\rm P}}}

\def\R{{\mathbb R}}

\def\Q{{\mathbb Q}}

\def\a{\alpha}

\def\ep{\varepsilon}
\def\eps{\varepsilon}

\def\E{{\mathbb E}}
\def\P{{\mathbb P}}

\makeatletter \@addtoreset{equation}{section} \makeatother
\newcommand {\qed}%
{%
    {}\hfill
    {}\hfill
    {$\square $}%
    \vspace {0.3cm}%
    \pagebreak [2]%
    \par
}%
\newenvironment{proof}[1]{%
    \vspace{0.3cm}%
    \pagebreak [2]%
    \par%
    \noindent {\bf  Proof~#1\ }}{\qed}%
\newenvironment{remark}{%
    \vspace{0.3cm} \pagebreak [2]%
    \par%
    \refstepcounter{proposition}
    \noindent%
   {\bf Remark~\theproposition\  }}{\ }%
\begin{document}
%

%\begin{titlepage}

\title{Sample Path Properties of Bifractional Brownian Motion}

\author{Ciprian A. Tudor \\ SAMOS-MATISSE, Centre d'Economie de La Sorbonne,\\Universit\'{e} de
Panth\'eon-Sorbonne Paris 1,\\90, rue de Tolbiac, 75634 Paris
Cedex 13, France.\\tudor@univ-paris1.fr\\
URL: \texttt{http://pageperso.aol.fr/cipriantudor/}\\
\vspace*{0.1in}\\
\and Yimin Xiao \thanks{Research partially supported by the NSF
grant DMS-0404729.}\\
Department of Statistics and Probability, \\
A-413 Wells Hall, Michigan State University,\\East Lansing, MI
48824,
USA.\\xiao@stt.msu.edu\\
URL: \texttt{http://www.stt.msu.edu/\~{}xiaoyimi}
 \vspace*{0.1in}} \maketitle

\begin{abstract}

Let  $B^{H, K}= \big\{B^{H, K}(t),\, t \in \R_+ \big\}$ be a
bifractional Brownian motion in $\R^d$. We prove that $B^{H, K}$
is strongly locally nondeterministic. Applying this property and a
stochastic integral representation of $B^{H, K}$, we establish
Chung's law of the iterated logarithm for $B^{H, K}$, as well as
sharp H\"older conditions and tail probability estimates for the
local times of $B^{H, K}$.

We also consider the existence and the regularity of the local
times of multiparameter bifractional Brownian motion
$B^{\overline{H}, \overline{K}}= \big\{B^{\overline{H},
\overline{K}}(t),\, t \in \R^N_+ \big\}$ in $\R^d$ using
Wiener-It\^o chaos expansion.\\
\end{abstract}

{\sc Running head:} Sample Path Properties of Bifractional
Brownian Motion\\

{\sc 2000 AMS Classification Numbers:} Primary 60G15, 60G17.\\


%\end{titlepage}
%\newpage
%\clearpage

\section{Introduction}

In recent years, there has been of considerable interest in
studying fractional Brownian motion due to its applications in
various scientific areas including telecommunications, turbulence,
image processing and finance. On the other hand, many authors have
proposed to use more general self-similar Gaussian processes and
random fields as stochastic models; see e.g. Addie et al. (1999),
Anh et al. (1999), Benassi et al. (2000), Mannersalo and Norros
(2002), Bonami and Estrade (2003), Cheridito (2004), Benson et al.
(2006). Such applications have raised many interesting theoretical
questions about self-similar Gaussian processes and fields in
general. However, contrast to the extensive studies on fractional
Brownian motion, there has been little systematic investigation on
other self-similar Gaussian processes. The main reasons for this,
in our opinion, are the complexity of dependence structures and
the non-availability of convenient stochastic integral
representations for self-similar Gaussian processes which do not
have stationary increments.

The objective of this paper is to fill this gap by developing
systematic ways to study sample path properties of self-similar
Gaussian processes. Our main tools are the Lamperti transformation
[which provides a powerful connection between self-similar
processes and stationary processes; see Lamperti (1962)] and the
strong local nondeterminism of Gaussian processes [see Xiao
(2005)]. In particular, for any self-similar Gaussian process $X
=\{X(t), t \in \R\}$, the Lamperti transformation leads to a
stochastic integral representation for $X$. We will show the
usefulness of such a representation in studying sample path
properties of $X$.

For concreteness, we only consider a rather special class of
self-similar Gaussian processes, namely, the bifractional Brownian
motions introduced by Houdr\'e and Villa (2003), to illustrate our
methods. %It will be clear that the arguments in this paper will be
%useful for studying self-similar Gaussian processes in general.
Given constants $H \in (0, 1)$ and $K \in (0, 1]$, the
bifractional Brownian motion (bi-fBm, in short) in $\R$ is a
centered Gaussian process $B^{H, K}_0 = \{B^{H, K}_0(t), t \in
\R_+\}$ with covariance function
\begin{equation}\label{Eq:covbifbm}
R^{H, K}(s, t) :=R(s, t) = \frac1 {2^K} \Big[ (t^{2H} + s^{2H})^K
- |t-s|^{2HK}\Big]
\end{equation}
and $B^{H, K}_0(0) = 0$.

Let $B^{H, K}_1, \ldots, B^{H, K}_d$ be independent copies of
$B^{H, K}_0$. We define the Gaussian process $B^{H, K} =
\big\{B^{H, K}(t),\, t \in \R_+ \big\}$ with values in $\R^d$ by
\begin{equation}\label{Eq:X}
B^{H, K}(t) = \Big(B^{H, K}_1(t), \ldots, B^{H, K}_d(t)\Big),
\qquad \forall t \in \R_+.
\end{equation}
By (\ref{Eq:covbifbm}) one can verify easily that $B^{H, K}$ is a
self-similar process with index $HK$, that is, for every constant
$a > 0$,
\begin{equation}\label{Eq:SS}
\Big\{ B^{H, K}(a t),\, t \in \R_+\Big\} \stackrel{d} {=} \Big\{
a^{HK} B^{H, K}(t), \, t \in \R_+\Big\},
\end{equation}
where $X \stackrel{d} {=} Y$ means the two processes have the same
finite dimensional distributions. Note that, when $K=1$,  $B^{H,
K}$ is the ordinary fractional Brownian motion in $\R^d$. However,
if $K \ne 1$, $B^{H, K}$ does not have stationary increments. In
fact, fractional Brownian motion is the only Gaussian self-similar
process with stationary increments [see Samorodnitsky and Taqqu
(1994)].

Russo and Tudor (2006) have established some properties on the
strong variations, local times and stochastic calculus of
real-valued bifractional Brownian motion. An interesting property
that deserves to be recalled is the fact that, when
$HK=\frac{1}{2}$, the quadratic variation of this process on $[0,
t]$ is equal to a constant times $t$. This is really remarkable
since as far as we know this is the only Gaussian self-similar
process with this quadratic variation besides Brownian motion.
Taking into account this property, it is natural to ask if the
bifractional Brownian motion $ B^{H, K}$ with $KH=\frac{1}{2}$
shares other properties with Brownian motion (from the sample path
regularity point of view). As it can be seen from the rest of the
paper, the answer is often positive: for example, the bi-fBm with
$HK=\frac{1}{2}$ and Brownian motion satisfy the same forms of
Chung's laws of the iterated logarithm and the H\"older conditions
for their local times.

The rest of this paper is organized as follows. In Section 2 we
apply the Lamperti transformation to prove the strong local
nondeterminism of $B_0^{H,K}$. This property will play essential
roles in proving most of our results. In Section 3 we derive small
ball probability estimates and a stochastic integral
representation for $B^{H,k}_0$. Applying these results, we prove a
Chung's law of the iterated logarithm for bifractional Brownian
motion.

Section 4 is devoted to the study of local times of one-parameter
bifractional Brownian motion and the corresponding $N$-parameter
fields. In general, there are mainly two methods in studying local
times of Gaussian processes: the Fourier analysis approach
introduced by Berman and the Malliavin calculus approach. It is
known that, the Fourier analysis approach combined with various
properties of {\em local nondeterminism} yields strong regularity
properties such as the joint continuity and sharp H\"older
conditions for the local times [see Berman (1973), Pitt (1978),
Geman and Horowitz (1980), Xiao (1997, 2005)]; while the Malliavin
calculus approach requires less conditions on the process and
establishes regularity of the local times in the sense of
Sobolev-Watanabe spaces [see Watanabe (1984), Imkeller et al.
(1995), Eddahbi et al. (2005)]. In this paper we make use of both
approaches to obtain more comprehensive results on local times of
bifractional Brownian motion and fields.

Throughout this paper, an unspecified positive and finite constant
will be denoted by $c$, which may not be the same in each
occurrence. More specific constants in Section $i$ are numbered as
$ c_{_{i, 1}}, c_{_{i, 2}}, \ldots$.
\vskip0.5cm

\section{Strong local nondeterminism}

The following proposition is essential in this paper. From its
proof, we see that the same conclusion holds for quite general
self-similar Gaussian processes.
\begin{proposition}\label{Prop:SLND}
For all constants  $0 < a < b$, $B_0^{H, K}$ is strongly locally
$\varphi$-nondeterministic on $I = [a, b]$ with $\varphi(r) = r^{2
HK}$. That is, there exist positive constants $c_{_{2,1}}$ and
$r_0 $ such that for all $t \in I$ and all $0 < r \le \min \{ t,\
r_0 \}$,
\begin{equation}\label{Eq:Cuzick}
{\rm Var} \Bigl(B_0^{H, K} (t) \big| B_0^{H, K}(s): s \in I,\ r
\le |s - t| \le r_0 \Bigr) \ge c_{_{2,1}}\, \varphi(r).
\end{equation}
\end{proposition}

\begin{proof}\  We consider the centered stationary Gaussian process
$Y_0 = \{Y_0(t), t \in \R\}$ defined through the Lamperti's
transformation [Lamperti (1962)]:
\begin{equation}\label{Eq:lam}
Y_0(t) = e^{-HK\,t}\, B_0^{H, K}(e^t),\quad \forall \, t \in \R.
\end{equation}
The covariance function $r(t) := \E\big(Y_0(0) Y_0(t)\big)$ is
given by
\begin{equation}\label{Eq:Ycov1}
\begin{split}
r(t)&= \frac 1 {2^K} e^{-HKt}\, \bigg[\big( e^{2Ht} +1\big)^K -
\big|e^t -
1\big|^{2HK}\bigg]\\
&= \frac 1 {2^K}\, e^{HKt}\, \bigg[\big( 1 + e^{-2Ht} \big)^K -
\big|1 - e^{-t} \big|^{2HK}\bigg].
\end{split}
\end{equation}
Hence $r(t)$ is an even function and, by (\ref{Eq:Ycov1}) and the
Taylor expansion, we verify that $r(t) = O( e^{-\beta t})$ as $t
\to \infty$, where $\beta = \min \{H(2-K), HK\}$. It follows that
$r(\cdot) \in L^1(\R)$. On the other hand,  by using
(\ref{Eq:Ycov1}) and the Taylor expansion again, we also have
\begin{equation}\label{Eq:Ycov2}
r(t)  \sim 1 - \frac 1 {2^K}\, |t|^{2HK}\qquad \hbox{ as } \ t \to
0.
\end{equation}

The stationary Gaussian process $Y_0$ is sometimes called the
Ornstein-Uhlenbeck process associated with $B_0^{H, K}$ [Note that
it does not coincide with the solution of the fractional Langevin
equation, see Cheridito et al. (2003) for a proof in the case
$K=1$].  By Bochner's theorem, $Y_0$ has the following stochastic
integral representation
\begin{equation}\label{Eq:YRep}
Y_0(t) = \int_{\R} e^{i \lambda t}\, W(d \lambda), \qquad
\forall\, t \in \R,
\end{equation}
where $W$ is a complex Gaussian measure with control measure
$\Delta$ whose Fourier transform is $r(\cdot)$. The measure
$\Delta$ is called the spectral measure of $Y$.

Since $r(\cdot) \in L^1(\R)$, so the spectral measure $\Delta$ of
$Y$ has a continuous density function $f(\lambda)$ which can be
represented as the inverse Fourier transform of $r(\cdot)$:
\begin{equation}\label{Eq:spden}
f(\lambda) = \frac 1 \pi \int_0^\infty r(t)\cos(t\lambda)\, d t.
\end{equation}
It follows from (\ref{Eq:Ycov2}), (\ref{Eq:spden}) and the
Tauberian theorem due to Pitman (1968, Theorem 5) [cf. Bingham et
al. (1987)] that
\begin{equation}\label{Eq:ftail}
f(\lambda) \sim c_{_{2,2}}\, |\lambda|^{-(1 + 2HK)} \quad \hbox{
as }\, \lambda \to \infty,
\end{equation}
where $c_{_{2,2}}>0$ is an explicit constant depending only on
$HK$. Hence, by a result of Cuzick and DuPreez (1982, Lemma 1)
[see also Xiao (2005) for more general results], $Y_0 = \{Y_0(t),
t \in \R\}$ is strongly locally $\varphi$-nondeterministic on any
interval $J = [-T, T]$ with $\varphi(r) = r^{2 HK}$ in the sense
that there exist positive constants $\delta$ and $c_{_{2,3}}$ such
that for all $t \in [-T, T]$ and all $r \in (0, |t| \wedge
\delta)$,
\begin{equation}\label{Eq:YLND}
{\rm Var} \bigl(Y_0(t) \big| Y_0(s): s \in J,\ r \le |s - t| \le
\delta \bigr) \ge c_{_{2,3}}\, \varphi(r).
\end{equation}

Now we prove the strong local nondeterminism of $B_0^{H, K}$ on
$I$. To this end, note that $B_0^{H, K}(t) = t^{HK} Y_0(\log t)$
for all $t > 0$. We choose $r_0 = a \delta$. Then for all $s, t
\in I$ with $r \le |s-t| \le r_0$ we have
\begin{equation}\label{Eq:log}
\frac r b \le \big|\log s - \log t\big| \le \delta.
\end{equation}
Hence it follows from (\ref{Eq:YLND}) and (\ref{Eq:log}) that for
all $t \in [a, b]$ and $r < r_0$,
\begin{equation}\label{Eq:BLND}
\begin{split}
& {\rm Var} \Bigl(B^{H, K}_0(t) \big|B^{H, K}_0 (s): s \in I,\ r
\le |s - t| \le r_0 \Bigr) \\
&= {\rm Var} \Bigl(t^{HK} Y_0(\log t) \big| s^{HK}\,Y_0(\log s): s
\in
I,\ r \le  |s - t| \le r_0 \Bigr) \\
&\ge t^{2HK}\, {\rm Var} \Bigl(Y_0(\log t) \big| Y_0(\log s): s
\in I,\ r \le  |s - t| \le r_0 \Bigr)\\
&\ge a^{2HK}{\rm Var} \Bigl(Y_0(\log t) \big| Y_0(\log s): s \in
I,\ r/b \le  |\log s - \log t| \le \delta \Bigr)\\
&\ge c_{_{2,4}}\, \varphi(r).
\end{split}
\end{equation}
This proves Proposition \ref{Prop:SLND}.
\end{proof}

For use in next section, we list two properties of the spectral
density $f(\lambda)$ of $Y$. They follow from (\ref{Eq:ftail}) or,
more generally, from (\ref{Eq:Ycov2}) and the truncation
inequalities in Lo\'eve (1977, p.209); see also Monrad and
Rootz\'en (1995).

\begin{lemma}\label{Lem:f}
There exist positive constants $ c_{_{2,5}}$ and $ c_{_{2,6}}$
such that for $u >1$,
\begin{equation}\label{Eq:f1}
\int_{|\lambda| < u} \lambda^2 f(\lambda)\, d \lambda \le
c_{_{2,5}}\, u^{2(1 -HK)}
\end{equation}
and
\begin{equation}\label{Eq:f2}
\int_{|\lambda| \ge u}  f(\lambda)\, d \lambda \le c_{_{2,6}}\,
u^{-2HK}.
\end{equation}
\end{lemma}

We will also need the following lemma from Houdr\'e and Villa
(2003).
\begin{lemma}\label{Lem:Var}
There exist positive constants $ c_{_{2,7}}$ and $ c_{_{2,8}}$ such
that for all $s, t \in \R _{+}$, we have
\begin{equation}\label{Eq:CompareVar}
 c_{_{2,7}}\, |t-s|^{2HK} \le \E\Bigl[\left( B^{H,K}_0(t )
 - B^{H,K}_0(s)\right)^2\Bigr] \le
c_{_{2,8}} \,|t - s|^{2HK}.
\end{equation}
\end{lemma}

\vspace{.02in}

\section{Chung's law of the iterated logarithm}
\label{Sec:LIL}

As applications of small ball probability estimates, Monrad and
Rootzen (1995), Xiao (1997) and Li and Shao (2001) established
Chung-type laws of the iterated logarithm for fractional Brownian
motion and other strongly locally nondeterministic Gaussian
processes with stationary increments. However, there have been no
results on Chung's LIL for self-similar Gaussian processes that do
not have stationary increments [Recall that the class of
self-similar Gaussian processes is large and fBm is the only such
process with stationary increments].

In this section, we prove the following Chung's law of the
iterated logarithm for bifractional Brownian motion in $\R$. It
will be clear that our argument is applicable to a large class of
self-similar Gaussian processes.

\begin{theorem}\label{Th:Chung}
Let $B^{H, K}_0 =\{B^{H, K}_0(t),\, t \in \R_+\}$ be a
bifractional Brownian motion in $\R$. Then there exists a positive
and finite constant $c_{_{3,1}}$ such that
\begin{equation}\label{Eq:LIL}
\liminf_{r \to 0} \frac{\max_{t \in [0, r]}\big|B^{H,K}_0(t)\big|}
{r^{HK}/(\log \log(1/r))^{HK}} = c_{_{3,1}} \qquad \hbox{ a.s.}
\end{equation}
\end{theorem}

In order to prove Theorem \ref{Th:Chung}, we need several
preliminary results. Lemma \ref{Th:SB} gives estimates on the
small ball probability of
$B^{H, K}_0$. %We refer to Li and Shao (2001) and Lifshits (1999)
%for extensive surveys on small ball probabilities, their
%applications and open problems.

\begin{lemma}\label{Th:SB}
There exist positive constants $c_{_{3,2}}$ and $c_{_{3,3}}$ such
that for all $t_0 \in [0, 1]$ and $ x \in (0, 1)$,
\begin{equation}\label{Eq:SB}
\exp\bigg(- \frac{c_{_{3,2}}} {x^{1/(HK)}} \bigg)\le
\P\left\{\max_{t \in [0, 1]} \big|B^{H,K}_0(t) -
B^{H,K}_0(t_0)\big| \le x \right\} \le \exp\bigg(-
\frac{c_{_{3,3}}} {x^{1/(HK)}} \bigg).
\end{equation}
\end{lemma}

\begin{proof}\ By Proposition \ref{Prop:SLND} and Lemma \ref{Lem:Var}, we see
that $B^{H, K}_0$ satisfies Conditions (C1) and (C2) in Xiao
(2005). Hence this lemma follows from Theorem 3.1 in Xiao (2005).
\end{proof}

Proposition \ref{Prop:01law} provides a zero-one law for ergodic
self-similar processes, which complements the results of Takashima
(1989). In order to state it, we need to recall some definitions.

Let  $X = \{X(t), t \in \R\}$ be a separable, self-similar process
with index $\kappa$. For any constant $a>0$, the scaling
transformation $S_{\kappa, a}$ of $X$ is defined by
\begin{equation}\label{Eq:Sa}
(S_{\kappa, a} X)(t) = a^{-\kappa} X(at), \qquad \forall t \in \R.
\end{equation}
Note that $X$ is $\kappa$-self-similar is equivalent to saying
that for every $a > 0$, the process $\{(S_{\kappa, a} X)(t), t \in
\R\}$ has the same finite dimensional distributions as those of
$X$. That is, for a $\kappa$-self-similar process $X$, a scaling
transformation $S_{\kappa, a}$ preserves the distribution of $X$,
and so the notion of ergodicity and mixing of $S_{\kappa, a}$ can
be defined in the usual way, cf. Cornfeld et al. (1982). Following
Takashima (1989), we say that a $\kappa$-self-similar process $X =
\{X(t), t \in \R\}$ is ergodic (or strong mixing) if for every $a
> 0, a \ne 1$, the scaling transformation $S_{\kappa, a}$ is
ergodic (or strong mixing, respectively). This, in turn, is
equivalent to saying that the shift transformations for the
corresponding stationary process $Y = \{Y(t), t \in \R\}$ defined
by $Y(t) = e^{-\kappa t} X(e^t)$ are ergodic (or strong mixing,
respectively).

\begin{proposition}\label{Prop:01law}
Let $X = \{X(t), t \in \R\}$ be a separable, self-similar process
with index $\kappa$. We assume that $X(0)=0$ and $X$ is ergodic.
Then for any increasing function $\psi: \R_+ \to \R_+$, we have
$\P(E_{\kappa, \psi}) = 0$ or 1, where
\begin{equation}\label{Eq:E}
E_{\kappa, \psi} = \left\{\omega: \hbox{ there exists $\delta > 0$
such that } \sup_{0\le s \le t} |X(s)| \ge t^{\kappa} \psi(t)\
\hbox{ for all } 0 < t \le \delta \right\}.
\end{equation}
\end{proposition}

\begin{proof}\ We will prove that for every $a > 0$, the event $E_{\kappa,
\psi}$ is invariant with respective to the transformation
$S_{\kappa, a}$. Then the conclusion follows from the ergodicity
of $X$.

Fix a constant $a > 0$ and $a \ne 1$. We consider two cases: (i)
$a > 1$ and (ii) $a < 1$. In the first case, since $\psi$ is
increasing, we have $\psi(au) \ge \psi(u)$ for all $u
> 0$. Assume that a.s. there is a $\delta > 0$ such that
\begin{equation}\label{Eq:case1}
\sup_{0\le s \le t} \left| X(s)\right|\ge t^{\kappa} \psi(t)\qquad
\hbox{ for all } 0 < t \le \delta,
\end{equation}
then
\begin{equation}\label{Eq:case12}
\sup_{0\le s \le t} \left| a^{-\kappa}\, X(as)\right| =
a^{-\kappa}\,\sup_{0\le s \le at} \left| X(s)\right| \ge
t^{\kappa} \psi(t)\quad  \hbox{ for all } 0 < t \le \delta/a.
\end{equation}
This implies that $E_{\kappa, \psi} \subset S^{-1}_{\kappa,
a}\big(E_{\kappa, \psi}\big)$. By the self-similarity of $X$,
these two events have the same probability, it follows that
$\P\big\{ E_{\kappa, \psi}\Delta S^{-1}_{\kappa, a}\big(E_{\kappa,
\psi}\big)\big\} = 0$. This proves that $E_{\kappa, \psi}$ is
$S_{\kappa, a}$-invariant and, hence, has probability 0 or 1.

In case (ii), we have $\psi(au) \le \psi(u)$ for all $u > 0$ and
the proof is similar to the above. If $S_{\kappa, a}X \in
E_{\kappa, \psi}$, then  we have $X \in E_{\kappa, \psi}$. This
implies $S^{-1}_{\kappa, a}\big(E_{\kappa, \psi}\big) \subset
E_{\kappa, \psi}$ and again $E_{\kappa, \psi}$ is $S_{\kappa,
a}$-invariant. This finishes the proof.
\end{proof}

By a result of Manuyama (1949) on ergodicity and mixing properties
of stationary Gaussian processes, we see that $B^{H, K}_0$ is
mixing. Hence we have the following corollary of Proposition
\ref{Prop:01law}.

\begin{corollary}\label{coro:01law}
There exists a constant $c_{_{3,4}} \in [0, \infty]$ such that
\begin{equation}\label{Eq:Chung01law}
\liminf_{t\to 0+} \frac{(\log \log 1/t)^{HK}} {t^{HK}}\, \max_{0
\le s\le t} \big|B^{H,K}_0(s)\big| = c_{_{3,4}},\qquad {\rm a.s.}
\end{equation}
\end{corollary}

\begin{proof}\ We take $\psi_c(t) = c\,\big(\log \log 1/t\big)^{-HK}$
and define $c_{_{3,4}} = \sup\big\{ c \ge 0:\,
\P\big\{E_{\kappa,\, \psi_c}\big\} = 1\big\}.$ It can be verified
that (\ref{Eq:Chung01law}) follows from Proposition
\ref{Prop:01law}.
\end{proof}

It follows from Corollary \ref{coro:01law} that Theorem
\ref{Th:Chung} will be established if we show $c_{_{3,4}} \in (0,
\infty)$. This is where Lemma \ref{Th:SB} and the following lemma
from Talagrand (1995) are needed. %When the function $\sigma^2(h)$ is well behaved [say,
%it is a power function], Lemma \ref{Lem:Tail} is well known.  See,
%for example, Xiao (2005) for a more general result.

\begin{lemma}\label{Lem:Tail}
Let $X = \{X(t), t \in \R\}$ be a centered Gaussian process in
$\R$ and let $S \subset \R$ be a closed set equipped with the
canonical metric defined by
$$ d(s, t) = \left[\E\big(X(s) - X(t)\big)^2\right]^{1/2}.$$
Then there exists a positive constants $c_{_{3, 5}}$ such that for
all $u > 0,$
\begin{equation}\label{Eq:Tail1}
\P\Biggl\{ \sup_{s, \ t \in S} |X(s) - X(t)| \ge c_{_{3, 5}}\,
\biggl(u + \int_0^D \sqrt{\log N_d(S, \ep)}\, d\ep \biggr)\Biggr\}
\le \exp\biggl( - \frac{u^2} {D^2}\biggr),
\end{equation}
where $N_d(S, \ep)$ denotes the smallest number of open $d$-balls
of radius $\ep$ needed to cover $S$ and where $D= \sup\{ d(s,t):\,
s,\, t \in S\}$ is the diameter of $S$.
\end{lemma}

Now we proceed to prove Theorem \ref{Th:Chung}.

\begin{proof} {\bf of Theorem \ref{Th:Chung}}\ We prove the lower bound first.
For any integer $n \ge 1$, let $r_n = e^{-n}$. Let $0 < \gamma <
c_{_{3,3}}$ be a constant and consider the event
\[
A_n = \bigg\{ \max_{0 \le s\le r_n} \big|B^{H,K}_0(s)\big| \le
\gamma^{HK} r_n^{HK}/(\log\log 1/r_n)^{HK}\bigg\}.
\]
Then the self-similarity of $B^{H,K}_0$ and Lemma \ref{Th:SB}
imply that
\begin{equation}\label{Eq:An1}
P\{A_n\} \le \exp\bigg(- \frac{c_{_{3,3}}} {\gamma} \log n \bigg)
= n^{- c_{_{3,3}}/\gamma}.
\end{equation}
Since $\sum_{n=1}^nP\{A_n\} < \infty$, the Borel-Cantelli lemma
implies
\begin{equation}\label{Eq:LIL-lb1}
\liminf_{n\to \infty} \frac{\max_{s \in [0,
r_n]}\big|B^{H,K}_0(s)\big|} {r_n^{HK}/(\log \log(1/r_n))^{HK}}
\ge c_{_{3,3}} \qquad \hbox{ a.s.}
\end{equation}
It follows from (\ref{Eq:LIL-lb1}) and a standard monotonicity
argument that
\begin{equation}\label{Eq:LIL-lb2}
\liminf_{r \to 0} \frac{\max_{t \in [0, r]}\big|B^{H,K}_0(t)\big|}
{r^{HK}/(\log \log(1/r))^{HK}} \ge c_{_{3,6}} \qquad \hbox{ a.s.}
\end{equation}

The upper bound is a little more difficult to prove due to the
dependence structure of $B^{H,K}_0$. In order to create
independence, we will make use of the following stochastic
integral representation of $B^{H, K}_0$: for every $t > 0$,
\begin{equation}\label{Eq:Rep1}
B^{H, K}_0(t) = t^{HK} \int_{\R} e^{i \lambda \log t}\,
W(d\lambda).
\end{equation}
This follows from the spectral representation (\ref{Eq:YRep}) of
$Y$ and its connection with $B^{H, K}_0$.

For every integer $n \ge 1$, we take
\begin{equation}\label{Eq:tn}
t_n = n^{-n} \ \ \ \hbox{ and }\ \ \ d_n = n^{\beta},
\end{equation}
where $\beta > 0$ is a constant whose value will be determined
later. It is sufficient to prove that there exists a finite
constant $ c_{_{3,7}}$ such that
\begin{equation}\label{Eq:UP}
\liminf_{n\to \infty} \frac{\max_{s \in [0,
t_n]}\big|B^{H,K}_0(s)\big|} {t_n^{HK}/(\log \log(1/t_n))^{HK}}
\le c_{_{3,7}} \qquad \hbox{ a.s.}
\end{equation}

Let us define two Gaussian processes $X_n$ and $\widetilde{X}_n$
by
\begin{equation}\label{Eq:Xn1}
X_n(t) = t^{HK} \int_{|\lambda| \in (d_{n-1}, d_n]} e^{i \lambda
\log t}\, W(d\lambda)
\end{equation}
and
\begin{equation}\label{Eq:Xn2}
\widetilde{X}_n(t) = t^{HK} \int_{|\lambda| \notin (d_{n-1}, d_n]}
e^{i \lambda \log t}\, W(d\lambda),
\end{equation}
respectively. Clearly $B^{H, K}_0(t) = X_n(t) + \widetilde
{X}_n(t)$ for all $t \ge 0$. It is important to note that the
Gaussian processes $X_n\, (n = 1, 2, \ldots)$ are independent and,
moreover, for every $n \ge 1,$ $X_n$ and $\widetilde{X}_n$ are
independent as well.

Denote $h(r) = r^{HK}\, \big(\log \log 1/r\big)^{-HK}$. We make
the following two claims:
\begin{itemize}
\item[(i).]\ There is a constant $\gamma > 0$ such that
\begin{equation}\label{Eq:UP1}
\sum_{n=1}^\infty \P\bigg\{\max_{s \in [0, t_n]} \big|X_n(s)\big|
\le \gamma^{HK}\, h(t_n)\bigg\} = \infty.
\end{equation}
\item[(ii).]\ For every $\eps> 0$,
\begin{equation}\label{Eq:UP2}
\sum_{n=1}^\infty \P\bigg\{\max_{s \in [0, t_n]}
\big|\widetilde{X}_n(s)\big| > \eps\, h(t_n)\bigg\} < \infty.
\end{equation}
\end{itemize}
Since the events in (\ref{Eq:UP1}) are independent, we see that
(\ref{Eq:UP}) follows from (\ref{Eq:UP1}), (\ref{Eq:UP2}) and a
standard Borel-Cantelli argument.

It remains to verify the claims (i) and (ii) above. By Lemma
\ref{Th:SB} and Anderson's inequality [see Anderson (1955)], we
have
\begin{equation}\label{Eq:UP3}
\begin{split}
\P\bigg\{\max_{s \in [0, t_n]} \big|X_n(s)\big| \le \gamma^{HK}\,
h(t_n)\bigg\}&\ge \P\bigg\{\max_{s \in [0, t_n]} \big|B^{H,
K}_0(s)\big| \le \gamma^{HK}\, h(t_n)\bigg\}\\
&\ge \exp\Big(- \frac{c_{_{3,2}}} {\gamma} \log (n \log n) \Big)\\
&= \big(n \log n\big)^{- c_{_{3,2}}/\gamma}.
\end{split}
\end{equation}
Hence (i) holds for $\gamma \ge c_{_{3,2}}$.

In order to  prove (ii), we divide $[0, t_n]$ into $p_n +1$
non-overlapping subintervals $J_{n, j} =[a_{n, j-1},\, a_{n, j}],$
$ (i = 0, 1, \ldots, p_n)$ and then apply Lemma \ref{Lem:Tail} to
$\widetilde{X}_n$ on each of $J_{n, j}$. Let $\beta> 0$ be the
constant in (\ref{Eq:tn}) and we take $J_{n, 0} = [0,\, t_n
n^{-\beta}]$. After $J_{n, j}$ has been defined, we take $a_{n,
j+1} = a_{n, j}(1 + n^{-\beta})$. It can be verified that the
number of such subintervals of $[0, t_n]$ satisfies the following
bound:
\begin{equation} \label{Eq:mJ}
p_n +1 \le c\, n^{\beta } \log n.
\end{equation}
Moreover, for every $j \ge 1$, if $s, t \in J_{n, j}$ and  $s <
t$, then we have $t/s - 1 \le n^{- \beta}$ and this yields
\begin{equation}\label{Eq:J1}
t-s \le s\, n^{-\beta}\  \ \ \hbox{ and }\ \ \log\Big(\frac t
s\Big) \le n^{-\beta}.
\end{equation}

Lemma \ref{Lem:Var} implies that the canonical metric $d$ for the
process $\widetilde{X}_n$ satisfies
\begin{equation}
d(s, t) \le c\, |s-t|^{HK}\quad \hbox{ for all }\ s, t > 0
\end{equation}
and $d(0, s) \le c\, t_n^{HK} \, n^{-\beta HK}$ for every $s \in
J_{n, 0}$. It follows that $D_0 := \sup\{ d(s,t); s,\ t \in J_{n,
0}\} \le c\, t_n^{HK} \, n^{-\beta HK}$ and
\begin{equation}\label{Eq:J01}
N_d(J_{n, 0}, \ep) \le c\, \frac{t_n\, n^{-\beta}} {\ep^{1/(HK)}}.
\end{equation}
Some simple calculation yields
\begin{equation}\label{Eq:J02}
\begin{split}
\int_0^{D_0} \sqrt{\log N_d(J_{n, 0}, \ep)}\ d\ep &\le
\int_0^{t_n^{HK} \, n^{-\beta HK}} \sqrt{\log
\Big(\frac{t_n \, n^{-\beta}} {\ep^{1/(HK)}}\Big)}\ d \ep\\
&= t_n^{HK} \, n^{-\beta HK} \, \int_0^1  \sqrt{\log
\Big(\frac{1} {u}\Big)}\ d u\\
&= c_{_{3,8}}\, t_n^{HK} \, n^{-\beta HK}.
\end{split}
\end{equation}
It follows from Lemma \ref{Lem:Tail} and (\ref{Eq:J02}) that
\begin{equation}\label{Eq:J03}
\P\bigg\{\max_{s \in J_{n, 0}} \big|\widetilde{X}_n(s)\big| >
\eps\, h(t_n)\bigg\} \le \exp\bigg(- c\, \frac{n^{2\beta HK}}
{\big(\log (n \log n)\big)^{2 HK}}\bigg).
\end{equation}

For every $1 \le j \le p_n$, we estimate the $d$-diameter of
$J_{n, j}$. It follows from (\ref{Eq:Xn2}) that for any $s, t \in
J_{n, j}$ with $s < t$,
\begin{equation}\label{Eq:Jj1}
\begin{split}
\E\Big(\widetilde{X}_n(s) - \widetilde{X}_n(t)\Big)^2 &=
\int_{|\lambda|\le d_{n-1}} \Big| t^{HK} \, e^{i \lambda \log t} -
s^{HK} \, e^{i \lambda \log s}\Big|^2 \,f(\lambda)\, d\lambda\\
&\qquad + \int_{|\lambda| > d_{n}} \Big| t^{HK} \, e^{i \lambda
\log t} - s^{HK} \, e^{i \lambda \log s}\Big|^2 \,f(\lambda)\,
d\lambda\\
&:= \EuScript{I}_1 + \EuScript{I}_2.
\end{split}
\end{equation}
The second term is easy to estimate: for all $s, t \in J_{n, j}$,
\begin{equation}\label{Eq:I2}
\EuScript{I}_2 \le 4\, t_n^{2HK}\, \int _{|\lambda| >
d_{n}}\,f(\lambda)\, d\lambda\le c_{_{3,9}}\, t_n^{2HK} \,
n^{-2\beta HK},
\end{equation}
where the last inequality follows from (\ref{Eq:f2}).

For the first term, we use the elementary inequality $1 - \cos x
\le x^2$ to derive that for all $s, t \in J_{n, j}$ with $s < t$,
\begin{equation}\label{Eq:I1}
\begin{split}
\EuScript{I}_1 &=\int_{|\lambda|\le d_{n-1}} \bigg[\big(t^{HK} -
s^{HK}\big)^2  + 2 t^{HK}\,s^{HK} \Big( 1 - \cos\big( \lambda \log
\frac t s\big)\Big)\bigg] \,f(\lambda)\, d\lambda\\
&\le s^{2HK}\, \Big(\frac t s - 1\Big)^{2HK} \int_{\R}
\,f(\lambda)\, d\lambda + 2 t^{2HK}\, \log^2\Big(\frac t s\Big)\,
\int_{|\lambda|\le d_{n-1}}  \lambda^2
\,f(\lambda)\, d\lambda\\
&\le c_{_{3,10}}\, t_n^{2HK} \, n^{- 2\beta HK},
\end{split}
\end{equation}
where, in deriving the last inequality, we have used (\ref{Eq:J1})
and (\ref{Eq:f1}), respectively.

It follows from (\ref{Eq:Jj1}), (\ref{Eq:I2}) and (\ref{Eq:I1})
that the $d$-diameter of $J_{n, j}$
satisfies
\begin{equation}\label{Eq:Jj2}
D_{j} \le c_{_{3,11}}\, t_n^{HK} \, n^{-\beta HK}.
\end{equation}
Hence, similar to (\ref{Eq:J03}), we use Lemma \ref{Lem:Tail} and
(\ref{Eq:Jj2}) to derive
\begin{equation}\label{Eq:Jj3}
\P\bigg\{\max_{s \in J_{n, j}} \big|\widetilde{X}_n(s)\big| >
\eps\, h(t_n)\bigg\} \le \exp\bigg(- c\, \frac{n^{2\beta HK}}
{\big(\log (n \log n)\big)^{2 HK}}\bigg).
\end{equation}

By combining (\ref{Eq:mJ}), (\ref{Eq:J03}) and (\ref{Eq:Jj3}) we
derive that for every $\ep > 0$,
\begin{equation}\label{Eq:UP4}
\begin{split}
\sum_{n=1}^\infty \P\Big\{\max_{s \in [0, t_n]}
\big|\widetilde{X}_n(s)\big| &> \eps\, h(t_n)\Big\} \le
\sum_{n=1}^\infty \sum_{j=0}^{p_n}\P\Big\{\max_{s \in J_{n, j}}
\big|\widetilde{X}_n(s)\big| > \eps\, h(t_n)\Big\}\\
&\le c\, \sum_{n=1}^\infty n^{\beta } \log n \, \exp\bigg(- c\,
\frac{n^{2\beta HK}} {\big(\log (n \log n)\big)^{2 HK}}\bigg)\\
& < \infty.
\end{split}
\end{equation}
This proves (\ref{Eq:UP2}) and hence the theorem.
\end{proof}

\begin{remark}\label{Coro:t0}
Let $t_0 \in [0, 1]$ be fixed and we consider the process $X =
\{X(t), t \in \R_+\}$ defined by $X(t)= B^{H,K}_0(t+t_0) -
B^{H,K}_0(t_0)$. By applying Lemma \ref{Th:SB} and modifying the
proof of Theorem \ref{Th:Chung}, one can show that
\begin{equation} \label{Eq:Chungt0}
c_{_{3,12}}^{-1} \le \liminf_{r \to 0} \frac{\max_{t \in [0,
r]}\big|B^{H,K}_0(t+t_0) - B^{H,K}_0(t_0)\big|} {r^{HK}/(\log
\log(1/r))^{HK}} \le c_{_{3,12}} \qquad \hbox{ a.s.},
\end{equation}
where $c_{_{3,12}}> 1$ is a constant depending on $HK$ only.
\end{remark}

Corresponding to Lemma \ref{Th:SB}, we can also consider the small
ball probability of $B^{H, K}_0$ under the H\"older-type norm. For
$\a \in (0, 1)$ and any function $y \in C_0([0, 1])$, we consider
the $\a$-H\"older norm of $y$ defined by,
\begin{equation}\label{Eq:fnorm}
\|y\|_\a = \sup_{s, t \in [0,1], s\ne t}\, \frac{|y(s)-y(t)|}
{|s-t|^\a}.
\end{equation}

The following proposition extends the results of Stolz (1996) and
Theorem 2.1 of Kuelbs, Li and Shao (1995) to bifractional Brownian
motion.

\begin{proposition}\label{Th:Stolz}
Let $B^{H, K}_0$ be a bifractional Brownian motion in $\R$ and $\a
\in (0, HK)$. There exist positive constants $c_{_{3,13}}$ and
$c_{_{3,14}}$ such that for all $ \ep \in (0, 1)$,
\begin{equation}\label{Eq:Stolz}
\exp\Big(- c_{_{3,13}}\, \ep^{ - 1/(HK-\a)}\Big)\le
\P\Big\{\|B^{H, K}_0 \|_\a \le \ep \Big\} \le \exp\Big(-
c_{_{3,14}}\, \ep^{- 1/(HK-\a)} \Big).
\end{equation}
\end{proposition}
\begin{proof}\ It follows from Theorem 3.4 in Xiao (2005).
\end{proof}

\section{Local times of bifractional Brownian motion}

This section is devoted to the study of the local times of the
bi-fBm  both in the one-parameter and multi-parameter cases. As we
pointed out in the Introduction there are essentially two ways to
prove the existence and regularity properties of local times for
Gaussian processes: the first is related to the Fourier analysis
and the local nondeterminism property; the second is based on the
Malliavin calculus and Wiener-It\^o chaos expansion. We will apply
the Fourier analysis approach for the one-parameter case and the
Malliavin calculus approach for the multiparameter case.

\subsection{The one-parameter case}

Let $B^{H, K} = \{B^{H, K}(t), t \in \R_+\}$ be a bifractional
Brownian motion with indices $H$ and $K$ in $\R^d$. For any closed
interval $I \subset \R_+$ and for any  $x\in \R^{d}$, the local
time $L(x,I)$ of $B^{H,K}$ is defined as the density of the
occupation measure $\mu_I$ defined by
$$
\mu_I (A)=\int_{I}1\!\!1_{A}\big(B^{H,K} (s)\big)\, ds,\qquad A\in
\mathcal{B}(\R^{d}).
$$
It can be shown [cf. Geman and Horowitz (1980) Theorem 6.4] that
the following \emph{occupation density formula} hods: for every
Borel function $g(t, x) \ge 0$ on $I \times \R^d$,
\begin{equation}\label{Eq:OCF}
\int_{I} g\big(t, B^{H,K}(t)\big)\, dt = \int_{\R^d}\int_I g(t, x)
L(x, dt)\, dx.
\end{equation}

Lemma \ref{Lem:Var} and Theorem 21.9 in Geman and Horowitz (1980)
imply that if $1/(HK) > d$ then $B^{H, K}$ has a local time $L(x,
t) := L(x,\, [0, t])$, where $(x, t) \in \R^d \times [0, \infty)$.
In fact, more regularity properties of $L(x, t)$ can be derived
from Theorem 3.14 in Xiao (2005) which we summarize in the
following theorem. Besides interest in their own right, such
results are also useful in studying the fractal properties of the
sample paths of $B^{H, K}$.

\begin{theorem}\label{Th:LHolder}
Let $B^{H, K} = \{B^{H, K}(t), t \in \R\}$ be a bifractional
Brownian motion with indices $H$ and $K$ in $\R^d$. If $1/(HK) >
d$, then the following properties hold:
\begin{itemize}
\item[(i)]\ $B^{H, K}$ has a local time $L(x, t)$ that is jointly
continuous in $(x, t)$ almost surely.

\item[(ii)]\ {\rm [Local H\"older condition]} For every $B \in
{\cal B}(\R)$, let $L^*(B) = \sup_{x\in \R^d} L(x, B)$ be the
maximum local time. Then there exists a positive constant $c_{_{4,
1}}$ such that for all $t_0 \in \R_+$,
\begin{equation}\label{Eq:localH}
\limsup_{r \to 0} \frac{L^*(B(t_0,r))} {\varphi_1(r)} \le c_{_{4,
1}}\qquad a.s.
\end{equation}
Here and in the sequel, $B(t, r) = (t-r, t+r)$ and $\varphi_1(r) =
r^{1 - HKd}(\log \log 1/r)^{HKd}$.

\item[(iii)]\ {\rm [Uniform H\"older condition]} For every finite
interval $I \subseteq \R$, there exists a positive finite constant
$c_{_{4, 2}}$ such that
\begin{equation}\label{Eq:uniformH}
\limsup_{r \to 0}\ \sup_{t_0 \in I} \frac{L^*(B(t_0,r))}
{\varphi_2(r)} \le c_{_{4, 2}} \qquad a.s.,
\end{equation}
where $\varphi_2(r) = r^{1 - HKd}(\log 1/r)^{HKd}.$
\end{itemize}
\end{theorem}

\begin{proof}\
By Proposition \ref{Prop:SLND} and Lemma \ref{Lem:Var}, we see
that the conditions of Theorem 3.14 in Xiao (2005) are satisfied.
Hence the results follow.
\end{proof}

The following states that the local H\"older condition for the
maximum local time is sharp.
\begin{remark}
By the definition of local times, we have that for every interval
$ Q \subseteq \R_+$,
\begin{equation}\label{Eq:LTChung}
|Q| = \int_{\overline{B^{H,K}(Q)}} L(x, Q)\, dx  \le L^*(Q) \cdot
\Big(\max_{s, t \in Q}  \big| B^{H, K}(s) - B^{H,
K}(t)\big|\Big)^d.
\end{equation}
By taking $Q = B(t_0, r)$ in (\ref{Eq:LTChung}) and using
(\ref{Eq:Chungt0}) in Remark \ref{Coro:t0}, we derive the lower
bound in the following
\begin{equation}\label{Eq:localHlower}
c_{_{4,3}} \le \limsup_{r \to 0} \frac{L^*(B(t_0,\,r))}
{\varphi_1(r)} \le c_{_{4, 4}}\qquad a.s.,
\end{equation}
where $c_{_{4,3}} >0$ is a constant independent of $t_0$ and the
upper bound is given by (\ref{Eq:localH}). A similar lower bound
for (\ref{Eq:uniformH}) could also be established by using
(\ref{Eq:LTChung}), if one proves that for every interval $I
\subseteq \R_+$,
\begin{equation}\label{Eq:UniChung}
\liminf_{r \to 0}\ \inf_{t \in I}\ \max_{s \in B(t, r)} \frac{|
B^{H,K}(s) - B^{H,K}(t)|} { r^{Hk}/(\log 1/r)^{HK}} \le c_{_{4,
5}} \quad \ a. s.
\end{equation}
\end{remark}

Theorem \ref{Th:LHolder} can be applied to determine the Hausdorff
dimension and Hausdorff measure of the level set $Z_x= \{t \in
\R_+: B^{H, K}(t) = x\}$, where $x \in \R^d$. See Berman (1972),
Monrad and Pitt (1987) and Xiao (1997, 2005). In the following
theorem we prove a uniform Hausdorff dimension result for the
level sets of $B^{H, K}$.

\begin{theorem}\label{Th:unilevel}
If $1/(HK) > d$, then with probability one,
\begin{equation}\label{Eq:leveldim}
\dim Z_x = 1 - HK d \ \ \hbox{ for all }\ x \in \R^d,
\end{equation}
where $\dim $ denotes Hausdorff dimension.
\end{theorem}
\begin{proof}\
It follows from Theorem 3.19 in Xiao (2005) that with probability
one,
\begin{equation}\label{Eq:leveldim2}
\dim Z_x = 1 - HK d \ \ \hbox{ for all }\ x \in {\mathcal O},
\end{equation}
where ${\mathcal O}$ is the random open set defined by
\[
{\mathcal O}   = \bigcup_{s, t \in \Q; \, s < t} \Big\{x \in
\R^d:\ L(x, [s, t]) > 0 \Big\}.
\]
Hence it only remains to show  ${\mathcal O} = \R^d$ a.s. For this
purpose, we consider the stationary Gaussian process $Y=\{Y(t), t
\in \R\}$ defined by $Y(t) = e^{-HKt} B^{H, K}(e^t)$, using the
Lamperti transformation.

Note that the component processes of $Y$ are independent and, as
shown in the proof of Proposition \ref{Prop:SLND}, they are
strongly locally $\varphi$-nondeterministic with $\varphi(r) =
r^{2HK}$. It follows from Theorem 3.14 in Xiao (2005) that $Y$ has
a jointly continuous local time $L_Y(x, t)$, where $(x, t ) \in
\R^d\times\R$. From the proof of Proposition \ref{Prop:SLND}, it
can be verified that $Y$ satisfies the conditions of Theorem 2 in
Monrad and Pitt (1987), it follows that almost surely for every $y
\in \R^d$, there exists a finite interval $J \subset \R$ such that
$L_Y(y, J) > 0$.

On the other hand, by using the occupation density formula
(\ref{Eq:OCF}), we can verify that the local times of $B^{H,K}$
and $Y$ are related by the following equation: for all $x \in
\R^d$ and finite interval $I=[a,\,b] \subset [0, \infty)$,
\begin{equation}\label{Eq:LTConn}
L(x, I) = \int_{[\log a, \, \log b]} e^{(1 - HK)s}\,
L_Y(e^{-HKs}\, x,\, ds).
\end{equation}
Hence, there exists a.s. a finite interval $I$ such that $L(0,
I)>0$. The continuity of $L(x, I)$ implies the a.s. existence of
$\delta > 0$ such that $L(y, I)> 0$ for all $y \in \R^d$ with $|y|
\le \delta$. Observe that the scaling property of $B^{H,K}$
implies that for all constants $c > 0$, the scaled local time
$c^{-(1 - HKd)} L(x, ct)$ is a version of $L(c^{-HK}x, t)$. It
follows that a.s. for every $x \in \R^d$, $L(x, J) > 0$ for some
finite interval $J \subset [0, \infty)$.
\end{proof}

Since there is little knowledge on the explicit distribution of
$L(0, 1)$, it is of interest in estimating the tail probability
$\P\{L(0, 1) > x\}$ as $x \to \infty$. This problem has been
considered by Kasahara et al. (1999) for certain fractional
Brownian motion and by Xiao (2005) for a large class of  Gaussian
processes. Our next result is a consequence of Theorem 3.20 in
Xiao (2005).

\begin{theorem}\label{Th:tail}
Let $B^{H, K} = \{B^{H, K}(t),\, t \in \R\}$ be a bifractional
Brownian motion in $\R^d$ with indices $H$ and $K$. If $1/(HK) >
d$, then for $x > 0$ large enough,
\begin{equation}\label{Eq:tail}
- \log \P\big\{L(0, 1) > x \big\} \asymp x^{HK},
\end{equation}
where $a(x) \asymp b(x)$ means $a(x)/b(x)$ is bounded from below
and above for $x $ large enough.
\end{theorem}

\begin{proof}\
By Proposition \ref{Prop:SLND} and Lemma \ref{Lem:Var}, we see
that the conditions of Theorem 3.20 in Xiao (2005) are satisfied.
This proves (\ref{Eq:tail}).
\end{proof}

Let us also note that the existence of the jointly continuous
version of the local time and the self-similarity allow us to
prove the following renormalization result. The case $d=1$ has
been proved in Russo and Tudor (2006).
\begin{proposition}
If $1/(HK) > d$, then for any integrable  function
$F:\mathbb{R}^{d}\to \mathbb{R}$,
\begin{equation} \label{ren}
t^{HK-1} \int_{[0, t]} F\big(B^{H,K}(u)\big)\, du \stackrel{(d)}
{\longrightarrow } \widetilde{F}\, L(0, 1) \quad \hbox{ as }\
{t\to \infty },
\end{equation}
where $\widetilde{F}= \int_{\R^d} F(x)\, dx $.
\end{proposition}
\begin{proof}\
It holds that
\begin{equation}\label{Eq:renormal1}
\int _{[0, t]} F \big(B^{H,K}(u)\big)\, du = t\int_{[0,1]}
F\big(B^{H,K} (tv)\big)\, dv  \stackrel{d} {=} t\int_{[0,1]}
F\left( t^{HK} B^{H,K}(v)\right)\, dv.
\end{equation}
By using the occupation density formula, we derive
\begin{equation}\label{Eq:renormal}
\int _{[0, t]} F \big(B^{H,K}(u))\, du=  t \int_{\mathbb{R}^{d}} F
\big( t^{HK}x \big)\, L (x,1)\, dx = t^{1-HK }\,
\int_{\mathbb{R}^{d}} F(y)\, L(y t^{-HK}, 1)\, dy.
\end{equation}
Since $y \mapsto L(y, 1)$ is almost surely continuous and bounded,
the dominated convergence theorem implies that, as $t \to \infty$,
the last integral in (\ref{Eq:renormal}) tends to $\widetilde{F}\,
L(0, 1)$ almost surely. This and (\ref{Eq:renormal1}) yield
(\ref{ren}).
\end{proof}

\subsection{Oscillation of bifractional Brownian motion}

The oscillations of certain classes of stochastic processes,
especially Gaussian processes,  in the measure space $([0,1],
\lambda_1 )$, where $\lambda_1 $ is the Lebesgue measure in $\R$,
have been studied, among others, by Wschebor (1992) and Aza\"is
and Wschebor (1996). The following is an analogous result for
bifractional Brownian motion.

\begin{proposition}
Let $B^{H,K}$ be a bi-fBm in $\R$ with indices $H \in (0,1)$ and
$K\in (0,1]$. For every $t\in [0,1]$, let
\begin{equation*}
Z_{\varepsilon }(t)= \frac{B^{H,K}(t+\varepsilon)
 -B^{H,K}(t)}{\varepsilon^{HK}}.
\end{equation*}
Then the following statements hold:
\begin{description}
\item{(i) } For every integer $k\geq 1$, almost surely,
\begin{equation*}
\int_{0}^{1} \big(Z_{\varepsilon } (t)\big)^{k}\, dt \to
\E(\rho^{k}) \qquad \hbox{ as }\ \varepsilon \to 0,
\end{equation*}
where $\rho$ is a centered normal random variable with variance
$\sigma ^{2}= 2^{1-K}$.

\item{(ii) } For every interval $J\subset [0,1]$, almost surely, for
every $x\in \mathbb{R}$
\begin{equation*}
\lambda_1 \{ t\in J:\, Z_{\varepsilon} (t) \leq x\} \to \lambda_1
(J)\, \P (\rho \leq x) \qquad \hbox{ as }\ {\varepsilon \to 0}.
\end{equation*}
\end{description}
\end{proposition}

\begin{proof}\  Let us denote
\begin{equation*}
Y^{\varepsilon, k}=\int_{0}^{1} \left(  Z_{\varepsilon }
(t)\right)^{k}\, dt.
\end{equation*}
It is sufficient to prove that
\begin{equation}\label{bY}
{\rm Var} \big(Y^{\varepsilon , k} \big) \leq c(k)\, \varepsilon
^{\beta} \hskip0.5cm \hbox{ for some } c(k) \ \hbox{ and } \beta
>0.
\end{equation}
Then the conclusions (i) and (ii) will follow as in Aza\"is and
Wschebor (1996) by the means of a Borel-Cantelli argument.

Note that
\begin{equation*}
{\rm Var} \big(Y^{\varepsilon , k}\big)  = \int_{0}^{1}
\int_{0}^{1} {\rm Cov} \big( Z _{\varepsilon }(u)^{k},\,
Z_{\varepsilon } (u)^{k}\big) \, dudv.
\end{equation*}
We will make use of the fact that for a centered Gaussian vector
$(U,V)$,
\begin{equation*}
{\rm Cov} \big(U^{k},\,V^{k}\big) = \sum_{1\leq p\leq k} c(p,k)
\big[ {\rm Cov}(U,\, V) \big]^{p} \, \big[ {\rm Var}(U) {\rm
Var}(V) \big]^{k-p}.
\end{equation*}
Since the random variable $Z_{\varepsilon}$ has clearly bounded
variance [cf. Lemma \ref{Lem:Var}], it suffices to show that for
every $1 \le p\le k$,
\begin{equation}\label{Eq:Suff}
\int_{0}^{1}\int_{0}^{1} \big[ \E\big( Z_{\varepsilon }(u)
Z_{\varepsilon }(v) \big) \big]^{p}\, dudv \leq c_{_{4, 6}} \,
\varepsilon^{\beta}
\end{equation}
We can write
\begin{equation*}
\begin{split}
\int_{0}^{1}\int _{0}^{1} \big[ \E\big( Z_{\varepsilon }(u)
Z_{\varepsilon }(v) \big) \big]^{p}\, dv du  &= 2\int_{0}^{1}\int
_{0}^{u} \I_{(u-v<\varepsilon )} \big[ \E\big( Z_{\varepsilon }(u)
Z_{\varepsilon }(v)\big) \big]^{p}\, dv du\\
&\qquad +2 \int_{0}^{1}\int_{0}^{u} \I_{(u-v \ge \varepsilon )}
\big[
\E \big( Z_{\varepsilon }(u) Z_{\varepsilon }(v)\big) \big]^{p}\, dv du\\
&: =A+B.
\end{split}
\end{equation*}
Clearly $A \leq c\, \varepsilon$, hence it suffices to bound the
term $B$. Note that
\begin{equation*}
\E\big( Z_{\varepsilon }(u) Z_{\varepsilon }(v) \big) =
\frac{1}{\varepsilon ^{2HK}}\int_{u-\varepsilon }^{u} \int
_{v-\varepsilon }^{v} \frac{\partial^{2}R}{\partial a
\partial b}\, db da.
\end{equation*}
Since
\begin{equation*}
\frac{\partial ^{2}R}{\partial a \partial b }(a,b)
=\frac{2HK}{2^{K}}\left[  \left( a^{2H} + b^{2H} \right) ^{K-2}
a^{2H-1}b^{2H-1} -(2HK-1) \vert a-b\vert ^{2HK-2}\right],
\end{equation*}
we have
\begin{equation*}
\begin{split}
B &\leq  c(p,H,K) \int _{0}^{1} \int _{0}^{u-\varepsilon } \left[
\frac{1}{\varepsilon ^{2HK}}\int _{u-\varepsilon }^{u} \int
_{v-\varepsilon }^{v}\left( a^{2H} + b^{2H} \right)^{K-2}
a^{2H-1}b^{2H-1}\, dbda \right]^{p}\, dvdu \\
& \qquad \quad + c(p,H,K) \int _{0}^{1} \int _{0}^{u-\varepsilon }
\left[ \frac{1}{\varepsilon ^{2HK}}\int _{u-\varepsilon }^{u} \int
_{v-\varepsilon }^{v} \vert a -b\vert ^{2HK-2}  dbda \right] ^{p} dvdu \\
&:= B_{1}+ B_{2}.
\end{split}
\end{equation*}
The term $B_{2}$ can be treated as in the fBm case [see Aza\"is and
Wschebor (1996), Proposition 2.1] and we get $B_{2} \leq c\,
\varepsilon^{\beta }$ for some constant $\beta >0$. Finally, since
$a^{2HK}+ b^{2HK} \geq a^{HK}b^{HK}$, we can write
\begin{equation*}
\begin{split}
B_{1} &\leq c(p,H,K) \int_{0}^{1} \int _{0}^{u-\varepsilon }
\left[ \frac{1}{\varepsilon ^{2HK}}\int _{u-\varepsilon }^{u} \int
_{v-\varepsilon }^{v} a^{HK-1}b^{HK-1}\, dbda \right]^{p}\, dvdu\\
&=c(p,H,K) \int_{0}^{1} \int_{0}^{u-\varepsilon } \left(
\frac{u^{HK} -(u-\varepsilon )^{HK}}{\varepsilon ^{HK}}\right)
^{p} \left( \frac{v^{HK} -(v-\varepsilon )^{HK}}{\varepsilon
^{HK}}\right)^{p}\, dvdu\\
& \le c\, \Bigg[\int_0^1\left( \frac{u^{HK} -(u-\varepsilon
)^{HK}}{\varepsilon ^{HK}}\right) ^{p}\, dv du\Bigg]^2.
\end{split}
\end{equation*}
A change of variable shows that $B_1 \le c\,
\varepsilon^{2(1-HK)}$. Combining the above yields
(\ref{Eq:Suff}). Therefore, we have proved (\ref{bY}), and the
proposition.
\end{proof}

\vskip0.5cm

The above result can be extended to obtain the almost sure weak
approximation of the occupation measure of the bi-fBm $B^{H,K}$ by
means of normalized number of crossing of $B^{H,K}_{\varepsilon}$,
where $B^{H,K}_{\varepsilon}$ represents the convolution of
$B^{H,K}$ with an approximation of the identity $\Phi
_{\varepsilon }(t)= \frac{1}{\varepsilon } \Phi \left(
\frac{t}{\varepsilon }\right) $ with $\Phi =\I_{[-1, 0]}$. If $g$
is a real function defined on an interval $I$, then the number of
crossing of level $u$ is
$$N_{u}(g,I)= \#  \{ t\in I, g(t)=u\},$$
where $\#E$ denotes the cardinality of $E$.

\begin{proposition}
Almost surely for every continuous function $f$ and for every
bounded interval $I \subset \R_+$,
\begin{equation*}
\left( \frac{\pi }{2}\right)^{1/2}\varepsilon ^{1-HK} \int
_{-\infty }^{\infty } f(u) N_{u} (B^{H,K}_{\varepsilon}, I)\, du
\to \int_{-\infty }^{\infty } f(u) L(u, I) \, du \qquad \hbox{ as
}\ {\varepsilon \to 0}.
\end{equation*}
\end{proposition}
\begin{proof}\
The arguments in Aza\"is and Wschebor (1996), Section 5, apply.
Details are left to the reader.
\end{proof}

\subsection{The multi-parameter case}

For any given vectors $H=(H_{1}, \ldots, H_{N}) \in (0, 1)^N$ and
$K =(K_{1}, \ldots ,K_{N}) \in (0, 1]^N$, an $(N, d)$-bifractional
Brownian sheet $B^{H, K} = \{B^{H, K}(t), t \in \R_+^N\}$ is a
centered Gaussian random field in $\R^d$ with i.i.d. components
whose covariance functions are given by
\begin{equation}\label{Eq:Covfunc}
\E\Big(B^{H, K}_1(s) B^{H, K}_1(t)\Big) = \prod_{j=1}^N \frac 1
{2^{K_j}} \bigg[ \Big(s_j^{2 H_j} + t_j^{2H_j}\Big)^{K_j} - |t_j -
s_j|^{2 H_j K_j}\bigg].
\end{equation}
It follows from (\ref{Eq:Covfunc}) that, similar to an $(N,
d)$-fractional Brownian sheet [cf. Xiao and Zhang (2002), Ayache
and Xiao (2005)], $B^{H, K}$ is operator-self-similar. However, it
does not have convenient stochastic integral representations which
have played essential r\^oles in the studies of fractional
Brownian sheets. Nevertheless, we will prove that the sample path
properties of $B^{H, K}$ are very similar to those of fractional
Brownian sheets, and we can describe the anisotropic properties of
$B^{H, K}$ in terms of the vectors $H$ and $K$.

We start with the following useful lemma.

\begin{lemma}\label{Lem:IncreVar}
For any $\ep> 0$, there exist  positive and finite constants
$c_{_{4,7}}$ and  $c_{_{4,8}}$ such that for all $s, t \in [\ep,
1]^N$,
\begin{equation}\label{Eq:IncVar}
c_{_{4,7}}\, \sum_{j=1}^N |s_j - t_j|^{2H_jK_j} \le
\E\bigg[\Big(B^{H, K}_1(s) - B^{H, K}_1(t)\Big)^2\bigg] \le
c_{_{4,8}}\, \sum_{j=1}^N |s_j - t_j|^{2H_jK_j},
\end{equation}
and
\begin{equation}\label{Eq:IncVar1}
c_{_{4,7}}\, \sum_{j=1}^N |s_j - t_j|^{2H_jK_j} \le {\rm
detCov}\Big(B^{H, K}_1(s), B^{H, K}_1(t)\Big) \le c_{_{4,8}}\,
\sum_{j=1}^N |s_j - t_j|^{2H_jK_j}.
\end{equation}
Here and in the sequel, detCov denotes determinant of the
covariance matrix.
\end{lemma}

\begin{proof}\ We will make use of the following
easily verifiable fact: For any Gaussian random vector $(Z_1,
Z_2)$,
\begin{equation}\label{Eq:Det}
{\rm detCov} (Z_1, Z_2) = {\rm Var} (Z_1) {\rm Var}(Z_2 | Z_1),
\end{equation}
where ${\rm Var} (Z_1)$ and  ${\rm Var}(Z_2 | Z_1)$ denote the
variance of $Z_1$ and the conditional variance of $Z_2$, given
$Z_1$, respectively.

By (\ref{Eq:Det}) we see that for all $s, t \in [\ep, 1]^N$,
\begin{equation}\label{Eq:ncVar4}
\begin{split}
{\rm detCov}\Big(B^{H, K}_1(s), B^{H, K}_1(t)\Big) &=\E\Big[B^{H,
K}_1(s)^2\Big]\,{\rm Var} \Big( B^{H, K}_1(t)\big|
B^{H, K}_1(s)\Big)\\
&\le \E\Big[B^{H, K}_1(s)^2\Big]\, \E\bigg[\Big(B^{H, K}_1(s) -
B^{H, K}_1(t)\Big)^2\bigg].
\end{split}
\end{equation}
Since ${\rm Var}\big(B^{H, K}_1(s)\big)$ is bounded from above and
below by positive and finite constants, it is sufficient to prove
the upper bound in (\ref{Eq:IncVar}) and the lower bound in
(\ref{Eq:IncVar1}).

When $N=1$, Lemma \ref{Lem:Var}, Proposition \ref{Prop:SLND} and
(\ref{Eq:Det}) imply that both (\ref{Eq:IncVar}) and
(\ref{Eq:IncVar1}) hold. Next we show that, if the lemma holds for
any $B^{H, K}$ with at most $n$ parameters, then it holds for
$B^{H, K}$ with $n+1$ parameters.

We verify the upper bound in (\ref{Eq:IncVar}) first. For any $s,
t \in [\ep, 1]^{n+1}$, let $s' = (s_1, \ldots, s_n, t_{n+1})$.
Then we have
\begin{equation}\label{Eq:IncVar2}
\begin{split}
\E\bigg[\Big(B^{H, K}_1(s) - B^{H, K}_1(t)\Big)^2\bigg] &\le 2
\E\bigg[\Big(B^{H, K}_1(s) - B^{H, K}_1(s')\Big)^2\bigg]\\
& \qquad + 2 \E\bigg[\Big(B^{H, K}_1(s') - B^{H,
K}_1(t)\Big)^2\bigg].
\end{split}
\end{equation}
For the first term, we note that whenever $s_1, \ldots, s_n \in
[\ep, 1]$ are fixed, $B^{H, K}$ is a (rescaled) bifractional
Brownian motion in $s_{n=1}$. Hence Lemma \ref{Lem:Var} implies
the first term in the right-hand side of (\ref{Eq:IncVar2}) is
bounded by $c\, |t_n - s_n|^{2 H_{n+1} K_{n+1}}$, where the
constant $c$ is independent of $s_1, \ldots, s_n \in [\ep, 1]$. On
the other hand, when $t_{n+1} \in [\ep, 1]$ is fixed, $B^{H, K}$
is a (rescaled) $(N, d)$-bifractional Brownian sheet. Hence the
induction hypothesis implies the second term in the right-hand
side of (\ref{Eq:IncVar2}) is bounded by $c\, \sum_{j=1}^n |t_j -
s_j|^{2 H_{j} K_{j}}$. This and (\ref{Eq:IncVar2}) together prove
the upper bound in (\ref{Eq:IncVar}).

Suppose the lower bound in (\ref{Eq:IncVar1}) holds for any $B^{H,
K}$ with at most $n$ parameters. For $N = n+1$, we write ${\rm
detCov} \big(B^{H, K}_1(s),\,  B^{H, K}_1(t)\big)$ as
\begin{equation}\label{Eq:IncVar3}
\begin{split}
&\prod_{j=1}^{n+1} t_j^{2 H_j K_j} \, s_j^{2 H_j K_j} -
\prod_{j=1}^{n+1} \frac{1} {2^{2 K_j}} \Big[ \big(t_j^{2H_j}
+ s_j^{2 H_j}\big)^{K_j} - |t_j - s_j|^{2 H_j K_j}\Big]^2 \\
& = \prod_{j=2}^{n+1} t_j^{2 H_j K_j} \, s_j^{2 H_j K_j} \bigg\{
s_{1}^{2 H_{1} K_{1}} \, t_{1}^{2 H_{1} K_{1}} - \frac{1} {2^{2
K_{1}}} \Big[\big( t_{1}^{2H_{1}}
+ s_{1}^{2 H_{1}}\big)^{K_{1}} - |t_{1} - s_{1}|^{2 H_{1} K_{1}}\Big]^2 \bigg\}\\
& \qquad + \frac{1} {2^{2 K_{1}}}\bigg[\big(t_{1}^{2H_{1}}
+ s_{1}^{2 H_{1}}\big)^{K_{1}} - |t_{1} - s_{1}|^{2 H_{1} K_{1}}\bigg]^2 \\
& \qquad \times \bigg\{\prod_{j=2}^{n+1} t_j^{2 H_j K_j} \, s_j^{2
H_j K_j}  -  \prod_{j=2}^{n+1} \frac{1} {2^{2 K_j}} \Big[
\big(t_j^{2H_j} + s_j^{2 H_j}\big)^{K_j} - |t_j - s_j|^{2 H_j
K_j}\Big]^2\bigg\}\\
&\ge c\, \sum_{j=1}^{n+1} |s_j - t_j|^{2H_jK_j}
\end{split}
\end{equation}
where the last inequality follows from the induction hypothesis.
This proves the lower bound in (\ref{Eq:IncVar1}).
\end{proof}

Applying Lemma \ref{Lem:IncreVar}, we can prove that many results
in Xiao and Zhang (2002), Ayache and Xiao (2005) on sample path
properties of fractional Brownian sheet also hold for $B^{H, K}$
as well. Theorem \ref{Th:LTexist} is concerned with the existence
of local times of $B^{H, K}$.

\begin{theorem}\label{Th:LTexist}
Let $B^{H,K}=\left\{ B^{H,K}(t),\, t\in \mathbb{R}_{+}^{N}\right\}
$ be an $(N,d)$-bifractional Brownian sheet with parameters $H\in
(0, 1)^N$ and $K \in (0, 1]^N$. If $d< \sum _{j=1}^{N}
\frac{1}{H_{j}K_{j}}$ then for any $N$-dimensional closed interval
$I \subset (0, \infty)^N$, $B^{H,K}$ has a local time $L(x,I)$,
$x\in \mathbb{R}^{d}$. Moreover, the local time admits the
following $L^{2}$-representation
\begin{equation}
\label{L2} L(x,I)= (2\pi ) ^{-d} \int_{\mathbb{R}^{d}
}e^{-i\langle y,x\rangle }\int_{I} e^{i\langle y, B^{H,K}(s)
\rangle}\, dsdy, \hskip0.5cm x\in \mathbb{R}^{d}.
\end{equation}
\end{theorem}

\begin{remark}
Although the existence of local times can also be proved by using
the Malliavin calculus [see Proposition \ref{ch} below],  we
prefer to provide a Fourier analytic proof because: 1) we can
compare in this way the two methods and 2) the above theorem gives
in addition the representation (\ref{L2}).
\end{remark}

\begin{proof}\
Without loss of generality, we may assume that $I=[\varepsilon, 1]
^{N}$ where $\varepsilon >0$. Let $\lambda_N$ be the Lebesgue
measure on $I$. We denote by $\mu$ the image measure of
$\lambda_N$ under the mapping $t \mapsto B^{H,K}(t)$. Then the
Fourier transform of $\mu$ is
\begin{equation}
\widehat{\mu} (\xi) = \int_{I} e^{i \l \xi,\ B^{H,K}(t) \r}\, dt.
\end{equation}
It follows from Fubini's theorem and (\ref{Eq:IncVar}) that
\begin{equation}\label{Eq:squareint}
\begin{split}
\E \int_{\R^d} \big|\widehat{\mu} (\xi) \big|^2 \, d \xi & =
\int_I \int_I \int_{\R^d} \E \Big(e^{i \l \xi,\ B^{H,K}(s) -
B^{H,K}(t)\r } \Big) \, d \xi \, ds dt \\
& = c\, \int_I \int_I \frac1 {\big[\E \big(B^{H,K}_1(s) -
B^{H,K}_1 (t) \big)^2\big]^{d/2}}\, ds dt \\
&\le c\, \int_I\int_I \frac1 {\big[ \sum_{j=1}^N |s_j - t_j|^{2H_j
K_j} \big]^{d/2}} \, ds dt.
\end{split}
\end{equation}
The same argument in Xiao and Zhang (2002, p. 214) shows that the
last integral is finite whenever $d< \sum _{j=1}^{N}
\frac{1}{H_{j}K_{j}}$. Hence, in this case, $\widehat{\mu} \in
L^2(\R^d)$ a.s. and Theorem \ref{Th:LTexist} follows from the
Plancherel theorem.
\end{proof}

\begin{remark}
Recently, Ayache, Wu and Xiao (2006) have shown that fractional
Brownian sheets have jointly continuous local times based on the
``sectorial local nondeterminism''. It would be interesting to
prove that $B^{H,K}$ is sectorially locally nondeterministic and
to establish joint continuity and sharp H\"older conditions for
the local times of $B^{H,K}$.
\end{remark}

\vskip0.4cm

Now we consider the Hausdorff and packing dimensions of the image,
graph and level set of $B^{H,K}$. In order to state our theorems
conveniently, we assume
\begin{equation}\label{Eq:H}
0 < H_1K_1 \le \ldots \le H_N K_N < 1.
\end{equation}

We denote packing dimension by $\dimp$; see Falconer (1990) for
its definition and properties. The following theorems can be
proved by using Lemma \ref{Lem:IncreVar} and the same arguments as
in Ayache and Xiao (2005, Section 3). We leave the details to the
interested reader.

\begin{theorem} \label{Th:dim}
With probability 1,
\begin{equation}\label{Eq:rangedim}
\dim B^{H, K} \big([0, 1]^N\big) = \dimp B^{H, K} \big([0,
1]^N\big) = \min \bigg\{ d; \quad \sum_{j=1}^N \frac1 {H_j K_j}
\bigg\}
\end{equation}
and
\begin{equation}\label{Eq:graphdim}
\begin{split}
& \qquad \dim {\rm Gr}B^{H, K}\big([0, 1]^N\big) =  \dimp {\rm
Gr}B^{H, K}\big([0, 1]^N\big) \\
&  \qquad = \left\{ \begin{array}{ll}
 \sum_{j=1}^N \frac1 {H_jK_j}  &\hbox{ if }\ \sum_{j=1}^N \frac 1 {H_jK_j} \le d, \\
\sum_{j=1}^k \frac{H_k K_k} {H_jK_j} + N-k + (1 - H_k K_k) d
&\text { if } \ \sum_{j=1}^{k-1} \frac 1 {H_j K_j} \le d <
\sum_{j=1}^{k} \frac 1 {H_j K_j}\, ,
 \end{array}\right.
 \end{split}
\end{equation}
where $\sum_{j=1}^{0} \frac 1 {H_jK_j} := 0$.
\end{theorem}

\vskip0.5cm

\begin{theorem} \label{Th:dimlevel}
Let $L_x = \{t \in (0, \infty)^N:\, B^{H, K} (t) = x\}$ be the
level set of $B^{H, K}$. The following statements hold:
\begin{itemize}
\item[(i)]\ If $\sum_{j=1}^N \frac1 {H_j} < d$, then for every $x\in
\R^d$ we have $L_x = \emptyset$ a.s.

\item[(ii)]\ If $ \sum_{j=1}^N \frac1 {H_j} > d$, then for any $x\in \R^d$
and $0 < \ep < 1$, with positive probability
\begin{equation}\label{Eq:dimlevel}
\begin{split}
\dim \big(L_x \cap [\ep, 1]^N\big) &= \dimp \big(L_x \cap [\ep,
1]^N\big) \\
&= \min \bigg\{ \sum_{j=1}^k \frac{H_k} {H_j} + N-k
- H_k d,\ 1 \le k \le N \bigg\} \\
&= \sum_{j=1}^k \frac{H_k} {H_j} + N-k - H_k d,  \qquad \text { if
} \ \sum_{j=1}^{k-1} \frac 1 {H_j} \le d < \sum_{j=1}^{k} \frac 1
{H_j}.
\end{split}
\end{equation}
\end{itemize}
\end{theorem}

\vskip0.2cm

\subsection{A Malliavin calculus approach}

Using the Malliavin calculus approach, we can study the local
times of  more general bifractional Brownian sheets. Consider the
$(N\times d)$-matrices
$$\overline{H}=(\overline{H}_1,\dots,\overline{H}_d) \ \ \ \mbox{ and
}\ \ \ \overline{K}=(\overline{K}_1,\dots,\overline{K}_d),
$$
where for any $i=1, \ldots, d$
\begin{equation*}
\overline{H}_i =\left( H_{i,1}, \dots, H_{i, N}\right) \ \ \
\mbox{ and } \ \ \ \overline{K}_i =\left( K_{i,1}, \dots, K_{i,
N}\right)
\end{equation*}
with $H_{i,j}\in (0,1) $ and $K_{i,j}\in (0, 1]$ for every $i=1,
\dots, d$ and $j=1, \dots, N$.

\smallskip

We will say that the Gaussian field $B^{\overline{H}, \overline{K}
}$ is an $(N,d)$-bifractional Brownian sheet with indices
$\overline{H}$ and $\overline{K}$ if
\begin{equation*}
B^{\overline{H}, \overline{K} }(t) = \left(
B^{\overline{H}_{1}}({t}), \ldots , B^{\overline{H}_{d}}({t})
\right), \hskip0.5cm t\in [0,\infty)^{N}
\end{equation*}
and for every $i=1, \dots, d$, the process
$\{B^{\overline{H}_{i}}({t}),\, t\in \R_+^{N}\}$ is centered and
has  covariance function
\begin{equation*}
\E\left( B^{\overline{H}_{i} , \overline{K}_{i} } ({t})
B^{\overline{H}_{i} , \overline{K}_{i} } (s)\right) =
R^{\overline{H}_{i} , \overline{K}_{i} }(s,t) = \prod _{j=1}^{N}
R^{H_{i,j}, K_{i,j}}(s_{j}, t_{j}).
\end{equation*}

As in subsection 4.1, the local time $L(x,\, t)$ ( $t\in {\mathbb
R}_+^{N}$ and $x\in \R^{d}$) of $B^{\overline{H}, \overline{K}}$
is defined as the density of the occupation measure $\mu_t,$
defined by
$$\mu_t (A)=\int_{[0,t]} \I_{A} \big(B^{\overline{H}, \overline{K}}(s)\big)\, ds,\quad A\in
\mathcal{B}(\R^{d}). $$

Formally, we can write
$$
L(x,\, t)=\int_{[0,t]}\delta_{x} \big(B^{\overline{H},
\overline{K}}(s) \big)\, ds,
$$
where $\delta _{x}$ denotes the Dirac function and $\delta _{x}(B_{s}^{%
\overline{H}, \overline{K}})$ is therefore a distribution in the
Watanabe sense (see Watanabe (1984)).

We need some notation. For $x \in \R$, let $p_{\sigma }(x)$ be the
centered Gaussian kernel with variance $\sigma >0$. Consider also
the Gaussian kernel on $\R^{d}$ given by
$$
p_{\sigma }^{d}(x)=\prod_{i=1}^d p_{\sigma }(x_{i}),\qquad
x=(x_{1},\ldots ,x_{d}) \in \R^d.
$$
Denote by $\mathbf{H}_{n}(x)$ the $n$--th Hermite polynomial
defined by $\mathbf{H}_0 (x)=1$ and for $n\geq 1$,
$$
\mathbf{H}_{n}(x)=\frac{(-1)^n}{n!}\exp
\Big({\frac{x^{2}}{2}}\Big)\, \frac{d^{n}}{dx^{n}}\exp \Big(
{-\frac{x^{2}}{2}}\Big),\quad  \ x\in \R.
$$

We will make use of the following technical lemma.
\begin{lemma}
For any $H\in (0,1)$ and $K\in (0, 1]$, let us define the function
\begin{equation*}
Q_{H,K}(z) = \frac{ R^{H,K}(1,z)}{z^{HK}}, \hskip0.5cm z\in (0,1]
\end{equation*}
and $Q_{H,K}(0)=0$. Then the function $Q_{H,K}$ takes values in
$[0,1]$, $Q_{H,K}(1)=1$ and it is strictly increasing. Moreover,
there exists $\delta>0$ such that for all $z\in (1-\delta,\, 1),$
\begin{equation}
\label{ineqQ} \left( Q_{H,K}(z) \right) ^{n} \leq \exp \left(
-c(\delta , H,K) n(1-z)^{2H}\right).
\end{equation}
\end{lemma}
\begin{proof}\
Clearly, the Cauchy-Schwarz inequality implies $0\leq Q_{H,K} (z)
\leq 1.$ Let us prove that the function $Q_{H,K}$ is strictly
increasing. By computing the derivative  $Q'_{H,K}(z)$ and
multiplying this by $z^{HK+1}$, we observe that this is equivalent
to show
\begin{equation}
\label{b3} (1-z) ^{2HK-1}(1+z) -(1+z^{2H})^{K-1}(1-z^{2H})> 0 \ \
\hbox{ for all }\ z \in (0, 1).
\end{equation}

If $HK\le \frac{1}{2},$ since $(1+z^{2H})^{K-1} \leq 1+z$, the left
side in (\ref{b3}) can be minorized by $(1+z^{2H})^{K}\left( (1-z)
^{2HK-1}-1+ z^{2H}\right)$ and this is positive since $(1-z)
^{2HK-1}\geq 1$.

If $HK>\frac{1}{2}$, we note that
\begin{eqnarray*}
&& (1-z) ^{2HK-1}(1+z) +(1+z^{2H})^{K-1}z^{2H}
\geq (1-z)(1+z) +(1+z^{2H})^{K-1}z^{2} \\
&& \geq (1+z^{2H})^{K-1}(1-z^{2}) +(1+z^{2H})^{K-1}z^{2}\geq
(1+z^{2H})^{K-1}.
\end{eqnarray*}
and this implies (\ref{b3}). Concerning the inequality
(\ref{ineqQ}), we note that
\begin{equation*}
Q_{H,K}(z)^{n}= \exp \left( n \log Q_{H,K}(z) \right) \geq \exp
\left( -n(1-Q_{H,K}(z))\right).
\end{equation*}
Now by Taylor's formula
\begin{equation*}
(1+z^{2H})^{K}z^{-HK} \leq 2^{K} + c(H,K,\delta) (1-z)^{2}
\end{equation*}
and therefore
\begin{eqnarray*}
Q_{H,K}(z)&&\leq 1 +c(H,K,\delta)
(1-z)^{2}-\frac{1}{2^{K}}(1-z)^{2HK}\\
&&\leq 1+c(H,K,\delta)(1-z)^{2HK}\delta
^{2-2HK}-\frac{1}{2^{K}}(1-z)^{2HK}.
\end{eqnarray*}
The conclusion follows as in the proof of Lemma 2 in Eddahbi et
al. (2005), since
\begin{equation*}
1-Q_{H,K}(z)\geq \frac{1}{2^{K}} (1-z) ^{2HK}(1-c(H,K,\delta))
\end{equation*}
for any $z\in (1-\delta , 1)$ with  $\delta$ close to zero and
with $c(H,K,\delta)$ tending to zero as $\delta \to 0$.
\end{proof}

\vskip0.5cm

The following proposition gives a chaotic expansion of the local
time of the $(N, d)$-bifractional Brownian sheet. The stochastic
integral $I_{n}(h)$ appeared below is the  multiple Wiener-It\^o
integral of order $n$ of the function $h$ of $nN$ variables with
respect to an $(N, 1)$ bifractional Brownian motion with
parameters $H=(H_{1}, \ldots, H_{N})$ and $K=(K_{1}, \ldots ,
K_{N})$ . Recall that such integrals can be constructed in general
on a Gaussian space [see, for example, Major (1981), or Nualart
(1995)]. We will only need the following isometry formula:
\begin{equation}\label{iso}
\E\left( I_{n} (\I_{[0,t]} ^{\otimes n})
I_{m}(\I_{[0,s]}^{\otimes m})\right) = R^{H,K}(t,s)^{n}
\I_{(n=m)}=\prod _{j=1}^{N} \left( R^{H_{j},K_{j}}(t_{j},s_{j})
\right) ^{n} \I_{(n=m)}
\end{equation}
for all $s, t\in \mathbb{R}^{N}_{+}$.

\begin{proposition}\label{ch}
For any $x\in \mathbb{R}^{d}$ and $t\in (0,\, \infty)^{N}$, the
local times $L(x,\,t)$ admits the following chaotic expansion
\begin{equation} \label{chaos}
L(x,\, t)=\sum_{n_1,\ldots ,n_d\geq 0}\int_{[0,\,t]} \prod_{i=1}^d
\frac{p_{\underline{s}^{2\overline{H}_i}\overline{K}_{i}}(x_i)}
{\underline{s}^{n_i\overline{H}_{i}\overline{K}_{i}}}\,
\mathbf{H}_{n_i} \Big(\frac{x_i}
{\underline{s}^{\overline{H}_i}}\Big) I_{n_i}^i(\I_{[0,s]}(\cdot
)^{\otimes n_i}) \;ds,
\end{equation}
where  $\underline{s}=s_1\cdots s_N$ and
$\underline{s}^{\overline{H}_i \overline{K}_{i}}=\prod_{j=1}^N
s_j^{H_{i,j}K_{i,j}}$. The integrals $I_{n_i}^i$ denotes the
multiple It\^{o} stochastic integrals with respect to the
independent $N$-parameter bifractional Brownian motion
$B^{\overline{H}_{i},\overline{K}_{i}}$.

\smallskip

Moreover, if $\sum_{j=1}^{N}\frac{1}{H_{j}^{\ast }K_{j}^{\ast
}}>d,$ where $H_{j}^{\ast}=\max\{H_{i,j}:\, i=1,\dots, d\}$ and  $
K_{j}^{\ast}=\max\{K_{i,j}:\, i=1,\dots, d\}$,  then $L(x,\,t)$ is
a random variable in $L^{2}(\Omega)$.
\end{proposition}

\begin{proof}\
The chaotic expression (\ref{chaos}) can be obtained similarly as
in Eddahbi et al. (2005) or Russo and Tudor (2006). It is based on
the approximation of the Dirac delta function by Gaussian kernels
with variance converging to zero. Let us evaluate the
$L^{2}(\Omega)$ norm of $L(x,\,t)$. By the independence of
components and the isometry of multiple stochastic integrals, we
obtain
\begin{equation}\label{Eq:L2n}
\begin{split}
\Vert L(x,\,t) \Vert _{2} ^{2} = \sum_{m\geq 0} \sum
_{n_{1}+\cdots + n_{d}=m} \int_{[0,t]} du \int _{[0,t]} dv \prod
_{i=1}^{d} \beta _{n_{i}}(u) \beta _{n_{i}}(v)
R^{\overline{H}_{i}, \overline{K}_{i}}(u,v)^{n_{i}},
\end{split}
\end{equation}
where
\begin{equation*}
\beta_{n_{i}}(u) = \frac{p_{\underline{s}^{2\overline{H}_i}
\overline{K}_{i}}(x_i)}
{\underline{s}^{n_i\overline{H}_{i}\overline{K}_{i}}} \,
\mathbf{H}_{n_i} \bigg(\frac{x_i} {\underline{s}^{\overline{H}_{i}
\overline{K}_{i} }}\bigg)
\end{equation*}
By Propositions 3 and 6 in Imkeller et al. (1995) [see also Lemma
11 in Eddahbi et al. (1995)], we have the bound
\begin{equation}
\label{b1} \beta_{n_{i}}(u) \beta _{n_{i}}(v) \leq c_{_{4, 9}}\,
\frac{1}{ (n_{i}\vee 1 ) ^{ \frac{8\beta -1}{6}}} \, \frac{1}{
\underline{u}^{n_{i}\overline{H}_{i}\overline{K}_{i} }
\underline{v}^{n_{i}\overline{H}_{i}\overline{K}_{i} }}
\end{equation}
for any $\beta \in [\frac{1}{4}, \frac{1}{2})$. Using the
inequality (\ref{b1}), we derive from (\ref{Eq:L2n}) that $\Vert
L(x,\,t) \Vert _{2} ^{2}$ is at most
\begin{equation}
\begin{split}
& c\, \sum_{m\geq 0} \sum _{n_{1}+\cdots +n_{d}=m} \bigg(\prod
_{i=1}^{d} \frac{1}{ (n_{i}\vee 1 )^{ \frac{8\beta -1}{6}} }\bigg)
\, \int _{[0,t]} du \int _{[0,u]}dv \prod_{i=1}^{d}
\prod_{j=1}^{N}\frac{ R^{H_{i,j}, K_{i,j}}(u_{j},v_{j})^{n_{i}}}{
(u_{j}v_{j}
)^{n_{i}H_{i,j} K_{i,j}} }\\
&= c \,\sum_{m\geq 0} \sum _{n_{1}+\cdots +n_{d}=m} \bigg( \prod
_{i=1}^{d} \frac{1}{ (n_{i}\vee 1 )^{ \frac{8\beta -1}{6}} }\bigg)
\, \prod_{j=1}^{N} \int_{0}^{t_{j}}u_{j}du_{j} \int_{0}^{1} \bigg(
\prod _{i=1}^{d} Q_{H_{i,j}, K_{i,j}}(z)
^{n_{i}} \bigg) dz \\
&= c_{_{4, 10}}\, \underline{t}^{2} \sum_{m\geq 0} \sum
_{n_{1}+\cdots + n_{d}=m} \bigg( \prod _{i=1}^{d} \frac{1}{
(n_{i}\vee 1 ) ^{ \frac{8\beta -1}{6}} }\bigg) \prod_{j=1}^{N}
\int_{0}^{1} \bigg( \prod _{i=1}^{d} Q_{H_{i,j}, K_{i,j}}(z)
^{n_{i}} \bigg)\, dz,
\end{split}
\end{equation}
where we used the change of variables $u_{j}=u_{j}$ and
$v_{j}=z_{j}u_{j}$. Using the above lemma and as in the proof of
Lemma 2 in Eddahbi et al. (2005), we can prove the bound
\begin{equation}\label{b2}
\int_{0}^{1} \bigg(\prod _{i=1}^{d} Q_{H_{i,j}, K_{i,j}}(z)
^{n_{i}} \bigg) dz \leq c_{_{4, 11}}\, m^{-\frac{1}{ 2H_{j}^{\ast}
K_{j}^{\ast}}}.
\end{equation}
Here $c_{_{4, 11}}=c_{_{4, 11}}(\overline{H}, \overline{K})$ depends
on $\overline{H}, \overline{K}$. Finally, (\ref{b2}) implies that
\begin{equation}\label{Eq:438}
\begin{split}
\Vert L(x,\,t) \Vert _{2} ^{2}&\leq  c_{_{4, 12}}\sum _{m\geq
1}\bigg( \prod_{j=1}^{N}m^{-\frac{1}{ 2H_{j}^{\ast}
K_{j}^{\ast}}}\bigg) \sum _{n_{1}+\cdots +n_{d}=m} \bigg(
\prod_{i=1}^{d} \frac{1}{ (n_{i}\vee 1 ) ^{ \frac{8\beta
-1}{6}} }\bigg)\\
&\leq c_{_{4, 13}}\sum _{m\geq 1} m^{ - \sum_{j=1}^{N}\frac{1}{
2H_{j}^{\ast} K_{j}^{\ast}} +d( 1-\frac{8\beta -1}{6})-1},
\end{split}
\end{equation}
where $c_{_{4, 12}}$ and $c_{_{4, 13}}$ depend on $\overline{H},
\overline{K}$  and $\underbar{t}$ only. The last series in
(\ref{Eq:438}) converges if
\begin{equation}
\label{Eq:Cond}
\sum_{j=1}^{N}\frac{1}{ 2H_{j}^{\ast} K_{j}^{\ast}}
> d \bigg( 1-\frac{8\beta -1}{6} \bigg).
\end{equation}
To conclude, observe that by choosing $\beta $ close to
$\frac{1}{2}$, $\sum_{j=1}^{N}\frac{1}{H_{j}^{\ast }K_{j}^{\ast
}}>d$ implies the required condition (\ref{Eq:Cond}).
\end{proof}

\vskip0.5cm

We recall that a random variable $F=\sum_{n}I_{n}(f_{n})$ belongs
to the Watanabe space $\mathbb{D}^{\alpha ,2}$ if
\begin{equation*}
\Vert F\Vert ^{2}_{\alpha ,2}: = \sum _{n\geq 0} (1+m)^{\alpha }\,
\Vert I_{n}(f_{n})\Vert _{2}^{2}<\infty.
\end{equation*}

\begin{corollary}
For any $t \in (0, \infty)^{N}$ and $x\in \mathbb{R}^{d}$, the
local time $L(x,\,t)$ of the $(N, d)$-bifractional Brownian sheet
$B^{\overline{H}, \overline{K}}$ belongs to the Watanabe space
$\mathbb{D}^{\alpha , 2}$ for every $0 < \alpha < \sum_{j=1}^{N}
\frac{1}{2H^{\ast }_{j}K^{\ast }_{j}}-\frac{d}{2}$.
\end{corollary}
\begin{proof}\
This is a consequence of the proof of Proposition \ref{ch}. Using
the computation contained there, we obtain for any $ \beta \in
[\frac{1}{4}, \frac{1}{2})$,
\begin{equation*}
\Vert L(x,\,t) \Vert ^{2}_{\alpha , 2}\leq c_{{4, 14}}(\overline{H},
\overline{K}, d,t) \sum_{m\geq 1} (1+m) ^{\alpha} \, m^{
d(1-\frac{8\beta -1}{6})-1 - \sum _{j=1}^{N} \frac{1}{2H^{\ast
}_{j}K^{\ast }_{j}}}
\end{equation*}
which is convergent if $\alpha <\sum _{j=1}^{N} \frac{1}{2H^{\ast
}_{j}K^{\ast }_{j}}-d(1-\frac{8\beta -1}{6})-1 - \sum _{j=1}^{N}
\frac{1}{2H^{\ast }_{j}K^{\ast }_{j}}$. Choosing $\beta $ close to
$\frac{1}{2}$, we get the conclusion.
\end{proof}

\bigskip

{\bf Acknowledgment}\ \ This work was initiated while both authors
were attending the Second Conference on Self-similarity and
Applications held during June 20--24, 2005, at INSA Toulouse,
France. We thank the organizers, especially Professor Serge Cohen,
for their invitation and hospitality.

\bibliographystyle{plain}
\begin{small}

\end{small}

\end{document}